\documentclass[12pt]{article}
\usepackage{latexsym,amsmath,amsthm}
\usepackage{amssymb}
\usepackage{graphicx}
\usepackage{appendix}
\usepackage{tikz}
\usepackage[colorlinks=true, linkcolor=blue]{hyperref}
\usepackage{bbm}

\newtheorem{thm}{Theorem}[section]

\newtheorem{claim}[thm]{Claim}

\theoremstyle{definition}
\newtheorem{prob}[thm]{Problem}

\theoremstyle{definition}
\newtheorem{question}[thm]{Question}

\theoremstyle{definition}

\theoremstyle{definition}

\theoremstyle{definition}

\theoremstyle{definition}

\theoremstyle{remark}

\theoremstyle{remark}

\usepackage[margin=1.05in]{geometry}

\usepackage{url}
\usepackage{microtype}

\newcommand{\weakbestdec}{0.347}
\newcommand{\weakbestinc}{0.280}

\begin{document}

\title{On the independent set sequence of a tree}

\author{Abdul Basit\thanks{Department of Mathematics,
Iowa State University, Ames IA; abasit@iastate.edu}\\David Galvin\thanks{Department of Mathematics,
University of Notre Dame, Notre Dame IN; dgalvin1@nd.edu. 
Supported in part by the Simons Foundation.
}}

\maketitle
 
\begin{abstract}
Alavi, Malde, Schwenk and Erd\H{o}s asked whether the independent set sequence of every tree is unimodal. Here we make some observations about this question. We show that for the uniformly random (labelled) tree, asymptotically almost surely (a.a.s.) the initial approximately 49.5\% of the sequence is increasing while the terminal approximately 38.8\% is decreasing. Our approach uses the Matrix Tree Theorem, combined with computation. We also present a generalization of a result of Levit and Mandrescu, concerning the final one-third of the independent set sequence of a K\"onig-Egerv\'ary graph. 

\end{abstract}

\section{Introduction}

A sequence $(a_0, a_1, \ldots, a_m)$ of real numbers is {\em unimodal} if there is $k$ such that
$$
a_0 \leq a_1 \leq \cdots \leq a_k \geq a_{k+1} \geq \cdots a_{m-1} \geq a_m.
$$
Unimodality is ubiquitous in combinatorics and algebra, see e.g. the survey papers \cite{Breaden, Brenti4, Stanley2}.

It is well-known that the matching sequence of any finite graph (the sequence whose $k$th term is the number of matchings with $k$ edges in the graph) is unimodal; this follows from the seminal theorem of Heilmann and Lieb \cite{HeilmannLieb} that the generating polynomial of the matching sequence has all real roots. In contrast, the {\em independent set sequence} of a graph $G$ --- the sequence whose $k$th term $i_k=i_k(G)$ is the number of independent sets (sets of mutually non adjacent vertices) of size $k$ in $G$ --- is not in general unimodal. Alavi, Malde, Schwenk and Erd\H{o}s \cite{AlaviErdosMaldeSchwenk} showed, in fact, that it can be arbitrarily far from unimodal, in a precise sense (see also \cite{BallGalvinHyryWeingartner}).

There are families of graphs for which the independent set sequence is known to be unimodal --- for example, claw-free graphs (graphs without an induced $K_{1,3}$), as first shown by Hamidoune \cite{Hamidoune}. In 1987 Alavi, Malde, Schwenk and Erd\H{o}s \cite{AlaviErdosMaldeSchwenk} posed a question about another very basic family:
\begin{question} \label{AMSE-tree-Q}
Is the independent set sequence of every tree unimodal? And what about every forest?
\end{question}
There have been numerous partial results, mostly exhibiting families of trees with unimodal independent set sequences, see e.g. \cite{BahlsBaileyOlsen, Bencs, GalvinHilyard, LevitMandrescu2, LevitMandrescu3, MandrescuSpivak, WangZhu, Zhu, ZhuChen, ZhuWang, Zhu2}. The unimodality of the independent set sequence of all forests on at most 25 vertices has been verified computationally \cite{Radcliffe, YosefMizrachiKadrawi}, but the full question remains stubbornly open. The best general result to date is due to Levit and Mandrescu. A {\em K\"onig-Egerv\'ary} graph is one in which the number of vertices is $\alpha+\mu$, where $\alpha$ is the size of the largest independent set and $\mu$ is the size of the largest matching (measured by number of edges). Bipartite graphs, and so in particular trees and forests, are K\"onig-Egerv\'ary. Levit and Mandrescu \cite{LevitMandrescu} show:
\begin{thm} \label{thm-LMpartialuni}
For a K\"onig-Egerv\'ary graph $G$,
$$
i_{\lceil (2\alpha-1)/3 \rceil} \geq i_{\lceil (2\alpha-1)/3 +1 \rceil} \geq \cdots \geq i_{\alpha-1} \geq i_\alpha.
$$ 
\end{thm} 
So the (non-zero part of the) independent set sequence of a tree is weakly decreasing for its last one-third. 
Theorem \ref{thm-LMpartialuni} is easily seen to be tight: the graph consisting of $\alpha$ vertex disjoint edges (and no other vertices) has independent set sequence which is weakly decreasing from exactly $i_{\lceil (2\alpha-1)/3 \rceil}$ on.

In this note we make a number of observations around Question \ref{AMSE-tree-Q}, the first of which is a generalization of Theorem \ref{thm-LMpartialuni}, showing that the theorem is more about graphs with independent sets of size at least half the number of vertices than about K\"onig-Egerv\'ary graphs. 
\begin{thm} \label{thm-decreasing-tail}
Let $G$ be a graph (not necessarily a tree or a K\"onig-Egerv\'ary graph) with $n$ vertices and maximum independent set size $\alpha$. The sequence $(i_k)_{k=\ell}^\alpha$ is weakly decreasing, where
$$
\ell = \left\lceil \frac{\alpha(n - 1)}{\alpha+n}\right\rceil.
$$   
If $\kappa$ satisfies $\alpha  \geq \kappa n$ then
\begin{equation}
\label{upper-tail}
\ell \leq \left\lceil \frac{\alpha}{1+\kappa} - \frac{\kappa}{1+\kappa} \right\rceil.
\end{equation}
\end{thm}
(See Section \ref{sec-proofs-tail} for the proof). The second part of Theorem \ref{thm-decreasing-tail} follows quickly from the first: if $\alpha  \geq \kappa n$ then 
$$
\frac{\alpha(n - 1)}{\alpha+n} \leq \frac{\alpha(n - 1)}{(1+\kappa)n} \leq \frac{\alpha}{1+\kappa} - \frac{\alpha}{(1+\kappa)n} \leq \frac{\alpha}{1+\kappa} - \frac{\kappa}{1+\kappa}. 
$$

Every $n$-vertex graph satisfies $\mu \leq n/2$, so every K\"onig-Egerv\'ary graph satisfies $\alpha \geq n/2$ (the converse of this is not true; for example, $K_3$ together with two isolated vertices has $3=\alpha \geq 5/2 =n/2$ but is not K\"onig-Egerv\'ary). Thus, taking $\kappa=1/2$ in (\ref{upper-tail}) we recover Theorem \ref{thm-LMpartialuni}. 

Theorem \ref{thm-decreasing-tail} gives no new information on Question \ref{AMSE-tree-Q}, the status of the independent set sequence of all trees, because there are trees with $\alpha = \lceil n/2 \rceil$. But for trees with $\alpha$ larger than $n/2$, it gives a decreasing tail longer than one-third of the length of the sequence. 

One obvious place to exploit this is in the study of the independent set sequence of the random uniform tree. Our model here is to select ${\bf T}$ uniformly from among the $n^{n-2}$ labelled trees on vertex set $\{1,\ldots,n\}$, and to consider the sequence $(X_0, X_1, \ldots, X_n)$ where $X_k$ is the number of independent sets of size $k$ in ${\bf T}$. To gain some evidence in favor of an affirmative answer to Question \ref{AMSE-tree-Q}, it is natural to ask whether $(X_0, X_1, \ldots, X_n)$ is a.a.s. (asymptotically almost surely --- with probability tending to $1$ as $n$ tends to infinity) unimodal. 

This seemingly simple question turns out to be quite intricate. It is easy, via the Matrix Tree Theorem, to establish $E(X_k)=\binom{n}{k}\left(1-\frac{k}{n}\right)^{n-1}$ (this was probably first observed by Bedrosian \cite{Bedrosian}), so that the sequence $(E(X_0), E(X_1), \ldots, E(X_n))$ is unimodal. One might then try to establish that with high probability the $X_k$'s fall in disjoint intervals centered around the $E(X_k)$'s, leading to a.a.s. unimodality. Unfortunately the variance of $X_k$ (which can also be explicitly calculated via the matrix tree theorem) turns out to be very large, typically much larger than $E(X_k)^2$, precluding a straightforward application of the second moment method.

Nonetheless, Theorem \ref{thm-decreasing-tail} allows us to say something about the decreasing tail of the independent set sequence of almost all trees, beyond what is given by Theorem \ref{thm-LMpartialuni}. Pittel \cite{Pittel}, tightening an earlier result of Meir and Moon \cite{MeirMoon}, established that for any $f(n)=\omega(1)$, a.a.s.
\begin{equation} \label{meir-moon}
\alpha({\bf T}) \in \left(\rho n-f(n)\sqrt{n},\rho n+f(n)\sqrt{n}\right)
\end{equation}
where $\alpha({\bf T})$ is the size of the largest independent set in ${\bf T}$, and $\rho \approx 0.5671$ is the unique real satisfying $\rho e^\rho=1$.
So, by Theorem \ref{thm-decreasing-tail},
a.a.s. $(i_k)$ is decreasing for $k \in [\rho n/(1+\rho),\rho n]$, or the last approximately $36\%$ of $[0,\rho n]$. 
Here we improve this. 
\begin{thm}
\label{thm-random-tree-dec-tail}
Let ${\bf T}$ be a uniformly random labelled tree on $n$ vertices, and let $X_k$ be the number of independent sets of size $k$ in ${\bf T}$. A.a.s. the sequence $(X_\ell, X_{\ell+1}, \ldots, X_n)$ is weakly decreasing, where  
$\ell=\weakbestdec n$. 
\end{thm} 
So a.a.s. the (non-zero part of the) independent set sequence of the uniform labelled tree is weakly decreasing for its terminal approximately $38.8\%$. See Section \ref{sec-proofs-mttapp} for the proof of Theorem \ref{thm-random-tree-dec-tail}. With some further computation it is likely that we could improve Theorem \ref{thm-random-tree-dec-tail}, but an improvement to $\ell=0.346n$ is beyond the reach of our current methods. (Problem \ref{prob-improve-wingard} suggests a possible direction of improvement.) 

\medskip

Our second observation around Question \ref{AMSE-tree-Q} concerns the start of the independent set sequence. Again, we begin with a general statement:  
\begin{thm} \label{thm-MT}
Let $G$ be a graph in which every maximal (by inclusion) independent set has size at least $\lambda$. Then the initial portion $(i_0, i_1, \ldots, i_{\lceil \lambda/2 \rceil})$ of the independent set sequence of $G$ is weakly increasing. 
\end{thm}
This is a straightforward generalization (see Section \ref{sec-proofs-head-plus-misc} for the short proof) of a result of Michael and Traves \cite{MichaelTraves}, who showed that if every independent set in $G$ is contained in an independent set of size $\alpha$ ($G$ is {\em well-covered}) then $i_0 \leq i_1 \leq \cdots \leq i_{\lceil \alpha/2 \rceil}$.

To connect this to the independent set sequence of a tree, we show:
\begin{thm} \label{thm-increasing-head}
Let $T$ be a tree with $n$ vertices and maximum independent set size $\alpha$. Every maximal independent set in $T$ has size at least $\left\lceil\frac{n-\alpha+1}{2}\right\rceil$, and so the initial portion $(i_0, i_1, \ldots, i_\ell)$ of the independent set sequence of $T$ is weakly increasing, where
$$
\ell =  \left\lceil \left\lceil\frac{n-\alpha+1}{2}\right\rceil/2\right\rceil =  \left\lceil \frac{n-\alpha+1}{4}\right\rceil.
$$
\end{thm}
(See Section \ref{sec-proofs-head-plus-misc} for the proof.) For example, if we know that $\alpha = \lceil n/2\rceil$ (its smallest possible value) then we get that the independent set sequence is increasing up to about $n/8$ or $0.25\alpha$. On the other hand, if we know that $\alpha = n-1$ (its largest possible value) then we get no information from Theorem \ref{thm-increasing-head}.

Recalling (\ref{meir-moon}), from Theorem \ref{thm-increasing-head} we can immediately say that a.a.s. the (non-zero part of the) independent set sequence of the uniform labelled tree is increasing for its initial about $19\%$, or up to about $0.1n$. By modifying the idea that goes into the proof of Theorem \ref{thm-MT}, we can improve this substantially.
\begin{thm} \label{thm-random-tree-inc-head}
Let ${\bf T}$ be a uniformly random labelled tree on $n$ vertices, and let $X_k$ be the number of independent sets of size $k$ in ${\bf T}$. A.a.s. the sequence $(X_0, X_1, \ldots, X_\ell)$ is weakly increasing, where  
$\ell=\weakbestinc n.$
\end{thm} 
So a.a.s. the (non-zero part of the) independent set of the uniform labelled tree is weakly increasing for its initial $49.5\%$. See Section \ref{sec-proofs-mttapp} for the proof of Theorem \ref{thm-random-tree-inc-head}. With some further computation it is likely that we could improve Theorem \ref{thm-random-tree-inc-head}, but an improvement to $\ell=0.281n$ is beyond the reach of our current methods. 
Using quite different methods Heilman \cite{Heilman} has recently shown that the independent set sequence of ${\bf T}$ is a.a.s. weakly increasing up to $0.265n$, or for the initial about $46\%$ of its non-zero part. 

\medskip

Our proofs of Theorems \ref{thm-random-tree-dec-tail} and \ref{thm-random-tree-inc-head} also give some information about the value of another well-studied graph parameter for the random labelled tree, namely the independent domination number. This is defined to be the cardinality of the smallest independent set that is also a dominating set (every vertex outside the set is adjacent to something in the set). Equivalently, it is the cardinality of the smallest independent set that is maximal (by inclusion). For any tree $T$ with $n$ vertices and $\ell=\ell(T)$ leaves it is known \cite{Favaron, Lemanska} that $i(T)$, the independent domination number of $T$, satisfies
\begin{equation} \label{eq-inddom}
\frac{n+2-\ell}{3} \leq i(T) \leq \frac{n+\ell}{3}.
\end{equation}
Since $\ell({\bf T})$ is concentrated around $n/e$, (\ref{eq-inddom}) says that with high probability $i({\bf T})$ is between about $0.210n$ and $0.456n$.

We can improve the lower bound. In Section \ref{sec-analysis} we present and analyse the asymptotics of the quantity $f(n,k,t)$, the expected number of independent sets of size $k$ in ${\bf T}$ that have exactly $t$ extensions to an independent set of size $k+1$. At $t=0$ this is exactly the expected number of maximal (by inclusion) independent sets of size $k$. The analysis of Section \ref{sec-analysis} (details omitted) shows that $f(n,\kappa n, 0)=o(1)$ for all $\kappa \leq 0.307$, and so (by Markov's inequality) we deduce that a.a.s. $i({\bf T})$ is at least $0.307n$.  
\begin{prob} \label{prob-inddomT}
Determine the a.a.s. behavior of $i({\bf T})$.
\end{prob}

\medskip

We end the introduction with a few further remarks about generalizations of Question \ref{AMSE-tree-Q}. Recall that there is a sequence of ever-stronger (first is implied by second, et cetera, but no reverse implications) conditions on a sequence $(a_0, \ldots, a_m)$ of positive terms:
\begin{itemize}
\item {\bf Unimodality}: $a_0 \leq a_1 \leq \cdots \leq a_k \geq a_{k+1} \geq \cdots \geq a_m$.
\item {\bf Log-concavity}: $a_k^2 \geq a_{k-1}a_{k+1}$ for $k=1, \ldots, m-1$.
\item {\bf Ordered log-concavity}: 
$$
a_k^2 \geq \left(1+\frac{1}{k}\right)a_{k-1}a_{k+1}
$$ 
for $k=1, \ldots, m-1$. (We say ``ordered'' because ordered log-concavity corresponds to the sequence $(k!a_k)_{k=0}^n$ being log-concave, and when $a_k$ counts objects each consisting of $k$ unordered elements, $k!a_k$ counts the same objects when also an order is put on the elements).
\item {\bf Ultra log-concavity}: 
$$
a_k^2 \geq \left(1+\frac{1}{k}\right)\left(1+\frac{1}{m-k}\right)a_{k-1}a_{k+1}
$$
for $k=1, \ldots, m-1$ (corresponding to the sequence $(a_k/\binom{m}{k})_{k=0}^m$, or equivalently $(k!(m-k)!a_k)_{k=0}^m$ being log-concave).  
\item {\bf Real roots}: $\sum_{k=0}^m a_kx^k$ has all real roots. 
\end{itemize}

Chudnovsky and Seymour \cite{ChudnovskySeymour} showed that the independent set sequence of a claw-free graph satisfies not just unimodality but the real roots property; on the other hand, the independent set sequence of trees does not in general satisfy ultra-log concavity, as witnessed by the star on four vertices. It is plausible, however that there is an affirmative answer to the following question: 
\begin{question} \label{GRad-tree-Q}
Is the independent set sequence of every tree ordered log-concave?
\end{question}
Radcliffe \cite{Radcliffe} has verified that every tree on up to 25 vertices has ordered log-concave independent set sequence (see also \cite{YosefMizrachiKadrawi}, where Yosef, Mizrachi and Kadrawi verify log-concavity for trees on up to 20 vertices).  

One reason to think about ordered log-concavity is that it has a very nice reformulation. For a graph $G$ with maximum independent set size $\alpha$, let ${\mathcal I}$ and ${\mathcal I}_k$ be the set of all independent sets of $G$, and the set of independent sets of size $k$, respectively. For $I \in {\mathcal I}$, denote by $e(I)$ the number of extensions of $I$ to an independent set of size $|I|+1$ (or: $e(I)$ is the number of vertices in $G$ that are neither in $I$ nor adjacent to anything in $I$). Denote by $e_k$ the average number of extensions of an independent set of size $k$, that is
$$
e_k = \frac{\sum_{I \in {\mathcal I}_k} e(I)}{i_k}.
$$
\begin{claim} \label{clm-ordered-lc}
The sequence $(i_k)_{k=0}^\alpha$ is ordered log-concave if and only if the sequence $(e_k)_{k=0}^{\alpha-1}$ is weakly decreasing.  
\end{claim}
(See Section \ref{sec-proofs-head-plus-misc} for the quick proof.) So Question \ref{GRad-tree-Q} is equivalent to:
\begin{question} \label{extensions-Q}
For every tree, is the sequence $(e_k)_{k=0}^{\alpha-1}$ weakly decreasing?
\end{question}

\medskip

Before turning to proofs of Theorems \ref{thm-decreasing-tail}, \ref{thm-random-tree-dec-tail}, \ref{thm-MT}, \ref{thm-increasing-head} and \ref{thm-random-tree-inc-head} and Claim \ref{clm-ordered-lc} (in Section \ref{sec-proofs}), we make a remark concerning the difference between Question \ref{AMSE-tree-Q} for trees versus forests. If $G$ has components $G_1, \ldots, G_k$, and component $G_\ell$ has independent set sequence $i^\ell = (i^\ell_0, i^\ell_1, \ldots)$, then the independent set sequence of $G$ is the convolution of the sequences $i^\ell$ --- that is, it is the coefficient sequence of the polynomial $\prod_{\ell=1}^k \sum_{j \geq 0} i^\ell_j x^j $. It is not in general the case that the convolution of unimodal sequences is unimodal, which means that Question \ref{AMSE-tree-Q} for trees is distinct from Question \ref{AMSE-tree-Q} for forests.
On the other hand, it is the case that the convolution of log-concave sequences is log-concave \cite{DavenportPolya}, which means that to establish the log-concavity of the independent set sequence of an arbitrary forest, it is sufficient to do so for an arbitrary tree. We do not at the moment know whether the convolution of ordered log-concave sequences is ordered log-concave.     

\section{Proofs} \label{sec-proofs}

\subsection{Proof of Theorem \ref{thm-decreasing-tail}} \label{sec-proofs-tail}

The proof follows from two old results. First, a theorem of Fisher and Ryan \cite{FisherRyan}: 
\begin{thm} \label{thm-FisherRyan}
For any graph $G$ with maximum independent set size $\alpha$, we have
$$
\left(\frac{i_1}{\binom{\alpha}{1}}\right)^\frac{1}{1} \geq \left(\frac{i_2}{\binom{\alpha}{2}}\right)^\frac{1}{2} \geq \left(\frac{i_3}{\binom{\alpha}{3}}\right)^\frac{1}{3} \geq \cdots \geq \left(\frac{i_{\alpha-1}}{\binom{\alpha}{\alpha-1}}\right)^\frac{1}{\alpha-1} \geq \left(\frac{i_\alpha}{\binom{\alpha}{\alpha}}\right)^\frac{1}{\alpha}.
$$ 
\end{thm}

Second, a theorem of Zykov \cite{Zykov}:
\begin{thm} \label{thm-Zykov}
For any graph $G$ with $n$ vertices and with maximum independent set size $\alpha$, and any $0 \leq k \leq \alpha$, we have
$$
i_k \leq \binom{\alpha}{k}\left(\frac{n}{\alpha}\right)^k.
$$ 
\end{thm}
(This is a corollary of a more general result that among all graphs on $n$ vertices with maximum independent set size $\alpha$, the one which maximizes the number of independent sets of size $k$ for each $0 \leq k \leq \alpha$ is the balanced union of $\alpha$ cliques.)

\medskip

\noindent {\bf Proof (of Theorem \ref{thm-decreasing-tail})}: 
From Theorem \ref{thm-FisherRyan} we see that for each $k\leq \alpha-1$ we have
$$
\left(\frac{i_k}{\binom{\alpha}{k}}\right)^\frac{1}{k} \geq \left(\frac{i_{k+1}}{\binom{\alpha}{k+1}}\right)^\frac{1}{k+1},
$$ 
so that if $i_{k+1} > i_k$ then
$$
\left(\frac{i_{k+1}}{\binom{\alpha}{k}}\right)^\frac{1}{k} > \left(\frac{i_{k+1}}{\binom{\alpha}{k+1}}\right)^\frac{1}{k+1}
$$
or
$$
i_{k+1} > \frac{\binom{\alpha}{k}^{k+1}}{\binom{\alpha}{k+1}^k}=\binom{\alpha}{k+1}\left(\frac{k+1}{\alpha-k}\right)^{k+1}.
$$
Now Theorem \ref{thm-Zykov} says
$$
i_{k+1} \leq \binom{\alpha}{k+1}\left(\frac{n}{\alpha}\right)^{k+1}
$$
from which we deduce
$$
\frac{n}{\alpha} > \frac{k+1}{\alpha-k}. 
$$
In summary: $i_{k+1} > i_k$ forces $k < (\alpha n - \alpha)/(\alpha+n)$, which implies Theorem \ref{thm-decreasing-tail}. \qed

\subsection{Proofs of Theorems \ref{thm-MT} and \ref{thm-increasing-head}, and of Claim \ref{clm-ordered-lc}}
\label{sec-proofs-head-plus-misc}

The proofs of Theorems  \ref{thm-random-tree-dec-tail}, \ref{thm-MT} and \ref{thm-random-tree-inc-head}, and of Claim \ref{clm-ordered-lc}, all have an element in common, which we introduce now.

Given a graph $G$ with maximum independent size $\alpha$, for $0 \leq j \leq \alpha-1$ denote by $B_j$ the bipartite graph with classes ${\mathcal I}_j$ (the set of independent sets  of size $j$ in $G$) and ${\mathcal I}_{j+1}$, with an edge joining $I \in {\mathcal I}_j$, $J \in {\mathcal I}_{j+1}$ if and only if $I \subseteq J$.

$B_j$ has $(j+1)i_{j+1}$ edges, since each independent set of size $j+1$ is an extension of exactly $j+1$ independent sets of size $j$. It also has $\sum_{I \in {\mathcal I}_j} e(I)$ edges, where as before $e(I)$ is the number of extensions of $I$ to an independent set of size $|I|+1$. So we have the identity
\begin{equation} \label{counting-edges}
\sum_{I \in {\mathcal I}_j} e(I) = (j+1)i_{j+1}
\end{equation}
for $j=0, \ldots, \alpha-1$.

\medskip

\noindent {\bf Proof} (of Theorem \ref{thm-MT}): For $k \leq \lambda$, each $I \in {\mathcal I}_{k-1}$ has $e(I) \geq \lambda-(k-1)$, since each such $I$ is in at least one independent set of size $\lambda$. From (\ref{counting-edges}) it follows that $(\lambda-(k-1))i_{k-1} \leq ki_k$, so that if 
$k \leq \lceil \lambda/2\rceil$ then $i_{k-1} \leq i_k$. \qed

\medskip

\noindent {\bf Proof} (of Theorem \ref{thm-increasing-head}): Let $K$ be a maximal independent set in $T$, of size $|K|$. 

Each of the $n-|K|$ vertices of $T-K$ must have at least one edge to $K$, so the subgraph induced by $T-K$ is a forest with $n-|K|$ vertices and at most $|K|-1$ edges, and so at least $n-2|K|+1$ components. It follows that $T-K$, and hence $T$, has an independent set of size at least $n-2|K|+1$. The result follows from $n-2|K|+1 \leq \alpha$. \qed

\medskip

\noindent {\bf Proof} (of Claim \ref{clm-ordered-lc}): 
From (\ref{counting-edges}) we have
$$
e_j = \frac{(j+1)i_{j+1}}{i_j},
$$
so that monotonicity of $(e_k)_{k=0}^{\alpha-1}$ is equivalent to
$$
\frac{i_1}{i_0} \geq \frac{2i_2}{i_1} \geq \cdots \geq \frac{ki_k}{i_{k-1}} \geq \frac{(k+1)i_{k+1}}{i_k} \geq \cdots \geq \frac{\alpha i_\alpha}{i_{\alpha-1}},
$$
which is in turn equivalent to ordered log-concavity of $(i_k)_{k=0}^\alpha$. \qed

\subsection{Proofs of Theorems \ref{thm-random-tree-dec-tail} and \ref{thm-random-tree-inc-head}} \label{sec-proofs-mttapp}

Theorem \ref{thm-MT} hinged on the identity (\ref{counting-edges}), which allows us to deduce that if every independent set of size $k$ has more than $k$ extensions to an independent set of size $k+1$, then $i_k \leq i_{k+1}$. For Theorem \ref{thm-random-tree-inc-head} we modify this to: if all but a vanishing proportion of independent sets of size $k$ have more than $k$ extensions to an independent set of size $k+1$, then a.a.s. $i_k \leq i_{k+1}$. Theorem \ref{thm-random-tree-dec-tail} depends on a similar statement, that if all but a vanishing proportion of independent sets of size $k$ have fewer than $k$ extensions to an independent set of size $k+1$, then a.a.s. $i_k \geq i_{k+1}$. In Section \ref{sec-setup} below we make these ideas precise. 
In Section \ref{sec-analysis} we do the necessary analysis on the inequalities presented in Section \ref{sec-setup}. In Section \ref{sec-mtt} we use the Matrix Tree Theorem to establish (\ref{eq-enkt}), the key identity used throughout.

\subsubsection{Key claim} \label{sec-setup}

Let $e(n,k,t)$ denote the probability that, in a uniformly chosen labelled tree on $[n]:=\{1,\ldots,n\}$, a particular set of size $k$ is an independent set and has exactly $t$ extensions to an independent set of size $k+1$. In Section \ref{sec-mtt} we use the Matrix Tree Theorem (and inclusion-exclusion) to establish
\begin{equation} \label{eq-enkt}
e(n,k,t) = \binom{n-k}{t}\sum_{\ell=0}^{n-k-t} (-1)^\ell \binom{n-k-t}{\ell}\left(1-\frac{k}{n}\right)^{t+\ell-1}\left(1-\frac{(k+t+\ell)}{n}\right)^k.
\end{equation}
Let $g_1(n,k)$ denote the expected number of independent sets of size $k$ that have no more than $k+ 1$ extensions to an independent set of size $k+ 1$, and let $g_2(n,k)$ denote the expected number of independent sets of size $k$ that have $k+1$ or more extensions to an independent set of size $k+ 1$. By linearity of expectation we have 
\begin{eqnarray*}
g_1(n,k)=\binom{n}{k}\sum_{t=0}^{k+1} e(n,k,t) & \mbox{and} & g_2(n,k)=\binom{n}{k}\sum_{t=k+1}^{n-k}  e(n,k,t).
\end{eqnarray*}

\begin{claim} \label{clm-prob-part}
Suppose that $n$ and $k$ with $k+2\leq n$ satisfy
\begin{equation} \label{main-inq}
g_1(n,k) \leq \frac{\binom{n-k+1}{k}}{n^2\log n}.
\end{equation}
Then all but a proportion $1/(n\log n)$ of trees on $[n]$ satisfy $i_k \leq i_{k+1}$. And if 
\begin{equation} \label{main-inq2}
g_2(n,k) \leq \frac{\binom{n-k}{k+1}}{n^2\log n}
\end{equation}
then all but a proportion $1/(n\log n)$ of trees on $[n]$ satisfy $i_{k+1} \leq i_k$.
\end{claim}

The proof will require the following result, which was possibly first observed by Wingard \cite[Theorem 5.1]{Wingard}):
\begin{thm} \label{thm-wingard}
For any tree $T$ on $n$ vertices, and for any $0 \leq k \leq n$,
$$
i_k(T) \geq i_k(P_n)~\left(= \binom{n-k+1}{k}\right),
$$
where $P_n$ is the path on $n$ vertices.
\end{thm}
In other words, within the family of trees, the path minimizes the number of independent sets of any size. (See Problem \ref{prob-improve-wingard} for further discussion).

\medskip

\noindent {\bf Proof} (of Claim \ref{clm-prob-part}): By Markov's inequality, under (\ref{main-inq}) all but a proportion at most $1/(n\log n)$ of trees on $[n]$ have no more than $\binom{n-k+1}{k}/n$ independent sets of size $k$ with no more than $k+ 1$ extensions to an independent set of size $k+ 1$. In what follows we work inside in this set ${\mathcal T}_1$ of trees.

As in the proofs of Theorem \ref{thm-MT} and Claim \ref{clm-ordered-lc}, for $T \in {\mathcal T}_1$ consider the bipartite graph $B_k$ with classes ${\mathcal I}_k$ (the set of independent sets of size $k$ in $T$) and ${\mathcal I}_{k+1}$, with an edge joining $I \in {\mathcal I}_k$, $J \in {\mathcal I}_{k+1}$ if and only if $I \subseteq J$. Recalling (\ref{counting-edges}) we have
\begin{equation} \label{eq-bip}
\sum_{I \in {\mathcal I}_k} e(I) = (k+1)i_{k+1}
\end{equation}
where $e(I)$ denotes the number of extensions of $I$ to an independent set of size $k+ 1$.

Now lower bounding $e(I)$ by $0$ if $I$ has no more than $k+ 1$ extensions to an independent set of size $k+ 1$, and by $k+2$ otherwise, we get
$$
\sum_{I \in {\mathcal I}_k} e(I) \geq (k+2)\left(i_k-\frac{\binom{n-k+1}{k}}{n}\right).
$$ 
Inserting into (\ref{eq-bip}) and using $k+2 \leq n$ yields  
\begin{eqnarray*}
(k+1)(i_{k+1}-i_k) & \geq & i_k-(k+2)\left(\frac{\binom{n-k+1}{k}}{n}\right) \\
& \geq & i_k - \binom{n-k+1}{k}. 
\end{eqnarray*}
That $i_k \geq \binom{n-k+1}{k}$ (completing the proof of the first part of the claim) follows from Theorem \ref{thm-wingard}. 

The proof of the second part of the claim is similar. By Markov, under (\ref{main-inq2}) all but a proportion at most $1/(n\log n)$ of trees on $[n]$ have no more than $\binom{n-k}{k+1}/n$ independent sets of size $k$ with $k+1$ or more extensions to an independent set of size $k+ 1$. In what follows we work inside in this set ${\mathcal T}_2$ of trees.

For $T \in {\mathcal T}_2$ we again consider the bipartite graph $B_k$. Upper bounding $e(I)$ by $n$ if $I$ has $k+ 1$ or more extensions to an independent set of size $k+ 1$, and by $k$ otherwise, we get
$$
(k+1)i_{k+1} = \sum_{I \in {\mathcal I}_k} e(I) \leq ki_k + \binom{n-k}{k+1}
$$ 
or
$$
k(i_{k+1}-i_k) \leq \binom{n-k}{k+1} - i_{k+1} \leq 0,
$$
the last inequality following from Theorem \ref{thm-wingard} applied to independent sets of size $k+1$.
\qed 

\subsubsection{Analysis} \label{sec-analysis}

To complete the proofs of Theorems \ref{thm-random-tree-dec-tail} and \ref{thm-random-tree-inc-head}, it remains to verify (\ref{eq-enkt}) (which we do in Section \ref{sec-mtt}), and to show that for all sufficiently large $n$, for all $k \geq \weakbestdec n$ (\ref{main-inq2}) holds (so that, by a union bound, all but a proportion at most $1/\log n$ of trees on $[n]$ satisfy $i_k \geq i_{k+1}$ for all $k \geq \weakbestdec n$), while for all $k \leq \weakbestinc n$ (\ref{main-inq}) holds. This section furnishes those verifications.

It will be convenient to introduce $f(n,k,t)$, the expected number of independent sets of size $k$ that have exactly $t$ extensions to an independent set of size $k+1$; we have
\begin{eqnarray}
f(n,k,t) & = & \binom{n}{k}e(n,k,t) \nonumber \\
& = & \binom{n}{k}\binom{n-k}{t} \times \nonumber \\
& & ~~~~~~\sum_{\ell=0}^{n-k-t} (-1)^\ell \binom{n-k-t}{\ell}\left(1-\frac{k}{n}\right)^{t+\ell-1}\left(1-\frac{(k+t+\ell)}{n}\right)^k \nonumber \\
& = & a(n,k,t)\sum_{\ell=0}^{n-k-t}\binom{n-k-t}{\ell}\left(\frac{k}{n}-1\right)^{\ell}\left(1-\frac{\ell}{n-(k+t)}\right)^k \label{alt-sum}
\end{eqnarray}
where
$$
a(n,k,t)= \binom{n}{k}\binom{n-k}{t}\left(1-\frac{k}{n}\right)^{t-1}\left(1-\frac{(k+t)}{n}\right)^k.
$$
The sum in (\ref{alt-sum}) is of the form $\sum_{\ell=0}^N \binom{N}{\ell}(x-1)^\ell\left(1-\frac{\ell}{N}\right)^k$, with $x < 1$, and is not easy to directly asymptotically analyze, since it is an alternating sum. However, with a little manipulation we can turn the sum into one involving only positive terms. Indeed, we have  
\begin{eqnarray}
\sum_{\ell=0}^N \binom{N}{\ell}(x-1)^\ell\left(1-\frac{\ell}{N}\right)^k
& = & \frac{(x-1)^N}{N^k}\sum_{\ell=0}^N \binom{N}{\ell}\left(\frac{1}{x-1}\right)^\ell\ell^k \nonumber \\ 
& = & \frac{1}{N^k}\sum_{j=1}^k {k \brace j}(N)_j x^{N-j} \label{eq-using-B}
\end{eqnarray}
where in the second line we use symmetry of the binomial coefficients and in the last line we use the identity
\begin{equation} \label{id-Boyadzhiev}
\sum_{\ell=0}^N \binom{N}{\ell}\ell^kz^\ell = \sum_{j=0}^k {k \brace j}(N)_j(1+z)^{N-j}z^j
\end{equation}
(see, e.g. \cite[Proposition 2.5]{Boyadzhiev}) with $z=1/(x-1)$. Here ${a \brace b}$ is a Stirling number of the second kind and $(a)_b$ is a falling power. 

While not necessary for our argument, let us observe that (\ref{id-Boyadzhiev}) admits a combinatorial proof. The left-hand side evidently counts triples consisting of a subset $S$ of a set of size $N$, a word of length $k$ over alphabet $S$, and a coloring of the letters of $S$ from a palette of $z$ colors. This collection of triples could also be determined by first choosing the number $j \in \{1, \ldots, k\}$ of distinct letters that appear in the word; then choosing the blocks in the word in which the same letter appears (${k \brace j}$ options); then choosing the letters that appear in each of these blocks (let $T$ be this set of letters; note $|T|=j$), and the colors that each of those letters receive ($(N)_jz^j$ options); and finally choosing the remainder of the selected letters (i.e., the rest of $S$), and their colors ($\sum_{X \subseteq [N]\setminus T} z^{|T|} = \sum_{i=0}^{n-j} \binom{N-j}{i}z^i = (1+z)^{N-j}$ options). This leads to a count of $\sum_{j=0}^k {k \brace j}N^{\underline{j}}z^j(1+z)^{N-j}$ for the number of triples.

With $N=n-k-t$ and $x=k/n$, (\ref{eq-using-B}) yields 
\begin{equation} \label{f-sum}
f(n,k,t)= \frac{a(n,k,t)(k/n)^{n-k-t}}{(n-k-t)^k}\sum_{j=1}^k {k \brace j}(n-k-t)_j \left(\frac{n}{k}\right)^j.
\end{equation}
Recalling Claim \ref{clm-prob-part}, our goal is to find the largest $k_1=k_1(n)$ and smallest $k_2=k_2(n)$ such that for all sufficiently large $n$ we have 
\begin{equation} \label{task}
\begin{array}{ll}
\binom{n-k+1}{k}^{-1}\sum_{t=0}^{k+1} f(n,k,t) \leq \frac{1}{n^2\log n} & \mbox{for all $k \leq k_1$}, \\ \binom{n-k}{k+1}^{-1}\sum_{t=k+1}^{n-k} f(n,k,t) \leq \frac{1}{n^2\log n} & \mbox{for all $k \geq k_2$},
\end{array}
\end{equation}
from which it follows that the independent set sequence of the uniform random labelled tree on $n$ vertices is almost surely weakly increasing up to $k_1$ and weakly decreasing from $k_2$ on.  

We first give a heuristic analysis of the right-hand side of (\ref{f-sum}). Setting $\kappa=k/n$ and $\tau=t/n$, and ignoring polynomial factors of $n$ in the approximations below, we have
$$
\frac{a(n,k,t)(k/n)^{n-k-t}}{(n-k-t)^k} \approx \left(\frac{\exp_2\left\{\left(H(\kappa)+(1-\kappa)H\left(\frac{\tau}{1-\kappa}\right)\right)\right\}(1-\kappa)^\tau(1-\kappa-\tau)^\kappa \kappa^{1-\kappa-\tau}}{(1-\kappa-\tau)^\kappa n^\kappa}\right)^n.
$$
(Here we use $\binom{a}{b} \approx \exp_2\{aH(b/a)\}$, where $H$ is the binary entropy function.)

To estimate the sum in (\ref{f-sum}) we start with the standard identity
\begin{equation} \label{stir-id}
\sum_{m \geq 0} {m \brace i} \frac{z^m}{m!} = \frac{\left(e^z-1\right)^i}{i!}
\end{equation}
from which we deduce 
\begin{equation} \label{ex-coeff}
\sum_{j=1}^k {k \brace j}(n-k-t)_j \left(\frac{n}{k}\right)^j = k! \left[z^k\right]\left(1+ \left(\frac{n}{k}\right)\left(e^z-1\right)\right)^{n-k-t},
\end{equation}
where $[z^k]$ is the operation that extracts the coefficient of $z^k$ from a power series in variable $z$.

We now appeal to a result of Good \cite[Theorem 6.1]{Good} (see also \cite[Theorem 2]{Gardy}) concerning the asymptotics of a coefficient of a high power of $z$ in the power series expansion of a high power of a power series in $z$. 
\begin{thm} \label{thm-good}
Suppose that $f(z)=\sum_{k=0}^\infty f_kz^k$ is a power series with positive coefficients and with infinite radius of convergence. Suppose that $N=N(r)$ ($r$ a natural number) is such that $N/r$ is bounded away from $0$ and from  infinity as $r \rightarrow \infty$. Then the implicit equation
\begin{equation} \label{rho-imp}
\frac{\rho f'(\rho)}{f(\rho)} = \frac{N}{r}
\end{equation}
defines a unique positive real $\rho=\rho(r)$, and
$$
\left|\left[z^N\right]\left(f(z)\right)^r - \frac{f(\rho)^r}{\sigma \rho^N \sqrt{2\pi r}}\right| \leq \frac{g(\rho)}{r} 
$$
as $r \rightarrow \infty$, 
where $\sigma = \sigma(r) > 0$ is defined by $\sigma^2= \left(\frac{f''(\rho)}{f(\rho)}-\frac{f'(\rho)^2}{f(\rho)^2} + \frac{f'(\rho)}{\rho f(\rho)}\right)\rho^2$ and $g$ is a continuous function. 
\end{thm}
See \cite[p. 868]{Good} for an explicit description of $g$. As observed in \cite{Gardy}, $\rho f'(\rho)/f(\rho)$ and $\sigma$ have probabilistic interpretations, that will by useful for us later: $\rho f'(\rho)/f(\rho)$ is the expectation of the probability distribution $X$, supported on the natural numbers, given by $P(X=k) \propto f_k\rho^k$, while $\sigma^2/\rho^2$ is the variance of $X$.     

Taking $f(z)=1+(1/\kappa)(e^z-1)$, defining $\rho$ implicitly via $\rho f'(\rho)/f(\rho) = \kappa/(1-\kappa-\tau)$, and using $k! \approx (k/e)^k$ we get from (\ref{ex-coeff}) that
$$
\sum_{j=1}^k {k \brace j}(n-k-t)_j \left(\frac{n}{k}\right)^j \approx \left(\frac{n^\kappa \kappa^\kappa f(\rho)^{1-\kappa-\tau}}{e^\kappa \rho^\kappa}\right)^n.
$$
It follows that $f(n,k,t) \approx C(\kappa,\tau)^n$ where
$$
C(\kappa,\tau) = \frac{\exp_2\left\{\left(H(\kappa)+(1-\kappa)H\left(\frac{\tau}{1-\kappa}\right)\right)\right\}(1-\kappa)^\tau \kappa^{1-\tau} f(\rho)^{1-\kappa-\tau}}{e^\kappa \rho^\kappa},
$$
so that, using $\binom{n-k-1}{k}, \binom{n-k}{k+1}\approx \exp_2\left\{(1-\kappa)H(\kappa/(1-\kappa))\right\}$ and recalling (\ref{task}), we get that (\ref{main-inq}) holds as long as
\begin{equation} \label{inq-heur}
\frac{\sup_{\tau \in [0,\kappa]} C(\kappa,\tau)}{\exp_2\{(1-\kappa)H(\kappa/(1-\kappa))\}} < 1
\end{equation}
and (\ref{main-inq2}) holds as long as
\begin{equation} \label{inq-heur2}
\frac{\sup_{\tau \in [\kappa,1]} C(\kappa,\tau)}{\exp_2\{(1-\kappa)H(\kappa/(1-\kappa))\}} < 1.
\end{equation}
We can computationally verify that (\ref{inq-heur}) holds for $\kappa \leq 0.280$ (but not for $\kappa=0.281$), and that (\ref{inq-heur2}) holds for $\kappa \geq 0.347$ (but not for $\kappa=0.346$), heuristically justifying the comments after the statements of Theorems \ref{thm-random-tree-dec-tail} and \ref{thm-random-tree-inc-head}.
We can also check that $C(\kappa,0) <1$ for all $\kappa \leq 0.307$, justifying the comment concerning the quantity $i({\bf T})$ made just before the statement of Problem \ref{prob-inddomT}. 

\medskip

Before making this heuristic analysis rigorous, we note one obvious place where both Theorems \ref{thm-random-tree-dec-tail} and \ref{thm-random-tree-inc-head} might be improved. Looking at the proof of Claim \ref{clm-prob-part}, we see that that the $\binom{n-k+1}{k}$ and $\binom{n-k}{k+1}$ on the right-hand sides of (\ref{main-inq}) and (\ref{main-inq2}), and so the $\binom{n-k+1}{k}^{-1}$ and $\binom{n-k}{k+1}^{-1}$ on the left-hand side of (\ref{task}), come directly from Theorem \ref{thm-wingard} (the number of independent sets of size $k$, or $k+1$, in any tree on $n$ vertices is at least the number of size $k$, or $k+1$, in the path on $n$ vertices). If we could replace $\binom{n-k+1}{k}$ and $\binom{n-k}{k+1}$ with something larger, then we would have that (\ref{task}) holds for larger $k_1=k_1(n)$ and smaller $k_2=k_2(n)$. 

If we were working with all trees, such a replacement would not be possible, since the bound in Theorem \ref{thm-wingard} is tight for some trees (e.g., for paths). But we are working only with almost all trees, and so for our purposes it would be enough to have an a.a.s. lower bound on the number of independent sets of size $k$ in the uniform random labelled tree. We can find such a bound, but it is not much larger than the deterministic bound --- for $k=\kappa n$, when $\binom{n-k+1}{k}$ and $\binom{n-k}{k+1}$ grow exponentially with $n$ (both with base $\exp_2\left\{H(\kappa/(1-\kappa))\right\}$), it is only larger by a polynomial in $n$, and this does not lead to any improvement in our results. What is needed for an improvement is an answer to the following problem:      
\begin{prob} \label{prob-improve-wingard}
Find an a.a.s. lower bound on the number of independent sets of size $k$ in the uniform random labelled tree of $n$ vertices, that is substantially better than $\binom{n-k+1}{k}$.
\end{prob}
``Substantially better'' here means that when $k=\kappa n$ (and $\kappa < 1/2$, to avoid trivialities) the bound should grow exponentially in $n$, with a base that is larger than $\exp_2\left\{H(\kappa/(1-\kappa))\right\}$.  

\medskip

The remainder of this section is devoted to making our heuristic analysis rigorous. We start with Theorem \ref{thm-random-tree-inc-head}, by analyzing (\ref{f-sum}) for $k \leq \weakbestinc n$ and $t \leq k+1$ (the range of values relevant for that theorem). Note that in proving Theorem \ref{thm-random-tree-inc-head} we may assume $k \geq 0.1n$, since we already know from Theorem \ref{thm-increasing-head} that a.a.s. the independent set sequence of the random tree is increasing up to $0.108n$. Also, we initially assume only that $k \leq 0.49n$.  
Recall that our specific goal is to establish
\begin{equation} \label{eq-inc-goal}
\sum_{t=0}^{k+1} \frac{f(n,k,t)}{\binom{n-k+1}{k}}
\leq \frac{1}{n^2\log n} 
\end{equation}
for all large enough $n$ and $k \leq \weakbestinc n$, where
$$
f(n,k,t)= \frac{a(n,k,t)(k/n)^{n-k-t}}{(n-k-t)^k}\sum_{j=1}^k {k \brace j}(n-k-t)_j \left(\frac{n}{k}\right)^j
$$
and
$$
a(n,k,t)= \binom{n}{k}\binom{n-k}{t}\left(1-\frac{k}{n}\right)^{t-1}\left(1-\frac{(k+t)}{n}\right)^k.
$$

We break $[0,n]$ into finitely many equal intervals, and for each $k$ and $t$ we upper bound the various terms that comprise $f(n,k,t)$ (and lower bound $\binom{n-k+1}{k}$) in terms of the upper and lower endpoints of the intervals in which $k$ and $t$ lie. This reduces the verification of (\ref{eq-inc-goal}) to a finite computation.  

So, let $M$ be some large, fixed, positive integer. Let sufficiently large $n$ be given (the need for $n$ to be large will arise around (\ref{inq-largen})). Let $k$ and $t$ be given, with $0.1n \leq k \leq 0.49n$ and $0 \leq t \leq k+1$. Let $1\leq p \leq M$ be that integer such that $(p-1)n/M \leq k < pn/M$ and let $1 \leq q \leq M$ be that integer such that $(q-1)n/M \leq t < qn/M$. We have the following straightforward bounds:
\begin{itemize}
    \item $\binom{n}{k} \leq \exp_2\left\{nH\left(\frac{k}{n}\right)\right\} \leq \exp_2\left\{nH\left(\frac{p}{M}\right)\right\}=A(p,M)^n$, where $H$ is the binary entropy function. Here, and throughout, we use the bound $\binom{a}{b} \leq \exp_2\left\{aH(b/a)\right\}$, and we also use that $p/M \leq 0.49$ (by our assumed upper bound $k \leq 0.49n$), so that $k/n \leq p/M$ implies $H(k/n) \leq H(p/M)$.
    \item $\binom{n-k}{t} \leq \exp_2\left\{(n-k)H\left(\frac{t}{n-k}\right)\right\} \leq \exp_2\left\{\left(1-\frac{(p-1)}{M}\right)nH\left(\frac{q}{M-p}\right)\right\}=B(p,q,M)^n$. Here we use that $q/(M-p) \leq 1/2$. 
    \item $\left(1-\frac{k}{n}\right)^{t-1} \leq \left(\frac{M}{M-p}\right)\left(1-\frac{(p-1)}{M}\right)^{\frac{(q-1)n}{M}}=\left(\frac{M}{M-p}\right)C(p,q,M)^n$. 
    \item $\left(1-\frac{(k+t)}{n}\right)^k \leq \left(1-\frac{(p+q-2)}{M}\right)^{\frac{(p-1)n}{M}}=D(p,q,M)^n$.
    \item $\left(\frac{k}{n}\right)^{n-k-t} \leq \left(\frac{p}{M}\right)^{n\left(1-\frac{(p+q)}{M}\right)}=E(p,q,M)^n$.
    \item $(n-k-t)^k \geq n^k\left(\left(1-\frac{(p+q)}{M}\right)^{\frac{pn}{M}}\right)=n^kF(p,q,M)^n$. We leave the $n^k$ term untouched here; it will be combined with a $k^k$ that will appear later. 
    \item Using the bound $\binom{a}{b} \leq \exp_2\{aH(b/a)\}/(a+1)$ (valid for all $0 \leq b \leq a$, $(a,b)\neq(0,0)$),
\begin{eqnarray*}
\binom{n-k+1}{k} & \geq & \frac{\exp_2\left\{(n-k+1)H\left(\frac{k}{n-k+1}\right)\right\}}{n+1} \\
& \geq & \frac{\exp_2\left\{n\left(1-\frac{p}{M}\right)H\left(\frac{p-1}{M-(p-2)}\right)\right\}}{n+1} = \frac{G(p,M)^n}{n+1}.
\end{eqnarray*}
\end{itemize}

Finally we deal with the sum in $f(n,k,t)$. Using (\ref{stir-id}) in the second line below
we have
\begin{eqnarray}
\sum_{j=1}^k {k \brace j}(n-k-t)_j \left(\frac{n}{k}\right)^j & \leq & \sum_{j=1}^k {k \brace j}(n-k-t)_j \left(\frac{M}{p-1}\right)^j \nonumber \\
& = & k! \left[z^k\right]\left(1+ \left(\frac{M}{p-1}\right)\left(e^z-1\right)\right)^{n-k-t} \nonumber \\
& \leq & \left(\frac{k}{e}\right)^k\left[z^k\right]\left(1+ \left(\frac{M}{p-1}\right)\left(e^z-1\right)\right)^{n-k-t} \nonumber \\
& \leq & \frac{k^k}{e^{\frac{(p-1)n}{M}}}\left[z^k\right]\left(1+ \left(\frac{M}{p-1}\right)\left(e^z-1\right)\right)^{n-k-t}. \label{need-good}
\end{eqnarray}
At this point we combine the factor of $n^k$ that appeared earlier (in the lower bound on $(n-k-t)^k$) with the factor $k^k/e^{(p-1)n/M}$ from (\ref{need-good}):
$$
\frac{k^k}{e^{\frac{(p-1)n}{M}}n^k} \leq \left(\frac{p}{eM}\right)^\frac{p-1}{M} = I(p,M)^n.
$$

We now use Theorem \ref{thm-good}, with $N=k$, $r=n-k-t$ and $f(z)=1+(M/(p-1))(e^z-1)$ (note that since $k \geq 0.1n$ we have $p>1$). For the $n, k, t$ we are considering we have $n-k-t \rightarrow \infty$, and we have that $k/(n-k-t)$ is confined to the constant interval
\begin{equation} \label{rho-int}
\left[\frac{p-1}{M-p-q+2},\frac{p}{M-p-q}\right].
\end{equation}
By the assumption $k \leq 0.49n$, and the fact that $t \leq k+1$, we have that the endpoints of the interval in (\ref{rho-int}) are positive, with the lower endpoint bounded away from $0$ and the upper endpoint bounded away from infinity.
Also note that $f(z)$ has power series about $0$ with all coefficients positive, and with infinite radius of convergence. So all hypotheses of Theorem \ref{thm-good} are satisfied. 

Following (\ref{rho-imp}), define $\rho=\rho(n,k,t) > 0$ implicitly by
\begin{equation} \label{imp-def}
\frac{k}{n-k-t} = \frac{\left(\frac{M}{p-1}\right)\rho e^\rho}{1+\left(\frac{M}{p-1}\right)\left(e^\rho-1\right)}.
\end{equation}
The conclusion of Theorem \ref{thm-good} is that
\begin{equation} \label{good-conclusion}
\left[z^k\right]\left(1+ \left(\frac{M}{p-1}\right)\left(e^z-1\right)\right)^{n-k-t} = \frac{\left(1+ \left(\frac{M}{p-1}\right)\left(e^\rho-1\right)\right)^{n-k-t}}{\sigma \rho^k \sqrt{2\pi(n-k-t)}}\left(1 \pm \frac{g(\rho)}{n-k-t}\right).
\end{equation}
Our goal now is to put an upper bound on the right-hand side of (\ref{good-conclusion}) in terms of $n, p$ and $q$, that has no dependence on $k$ or $t$ (and so also no dependence on $\rho$).

A crucial observation is that if $c > 0$ and $y=y(x)$ is the unique positive solution to $x=cye^y/(1+c(e^y-1))$, then for $x \in (0,\infty)$ it holds that $y$ is monotone increasing. Indeed, we have
$$
\frac{dy}{dx} = \frac{(1+c(e^y-1))^2}{ce^y(1+c(e^y-1)+y(1-c))},
$$
so we need only check the positivity of $1+c(e^y-1)+y(1-c))$ for positive $y$; at $y=0$ this expression takes value $1$, and its derivative is $1+c(e^y-1)$, which is positive for positive $y$.
It follows that
$$
0 < \rho_{\rm min} = \rho\left(n,\frac{(p-1)n}{M},\frac{(q-1)n}{M}\right) \leq \rho(n,k,t) \leq \rho\left(n,\frac{pn}{M},\frac{qn}{M}\right) = \rho_{\rm max},
$$
where notice that $\rho_{\rm min}$ and $\rho_{\rm max}$ depend only on $p$ and $q$ (and $M$), but not on $n$. Noting that $1+(M/(p-1))(e^x-1)\geq 1$ for all $x \geq 0$ it follows that
$$
\frac{\left(1+\left(\frac{M}{p-1}\right)(e^\rho-1)\right)^{n-k-t}}{\rho^k} \leq \left(\frac{\left(1+\left(\frac{M}{p-1}\right)(e^{\rho_{\rm max}}-1)\right)^{1-\frac{(p+q-2)}{M}}}{\min\{\rho_{\rm min}^{p/M}, \rho_{\rm min}^{(p-1)/M}\}}\right)^n =J(p,q,M)^n.
$$
To explain the denominator: observe that in putting an upper bound on $1/\rho^k$ we have to pay attention to whether $\rho$ is smaller than $1$ --- in which case we should use an upper bound for $k$ --- or $\rho>1$ --- in which case we should use a lower bound for $k$.

We now argue that $\sigma=\sigma(n,k,t)$ is bounded below by a positive constant depending only on $p, q$ and $M$. This follows from the fact that, as observed after the statement of Theorem \ref{thm-good}, $\sigma/\rho$ is the standard deviation of a probability distribution (supported on ${\mathbb N}$) that is not almost surely constant, so $\sigma >0$. We may lower bound $\sigma$ by the minimum value it attains as the left-hand side of (\ref{imp-def}) varies over the (\ref{rho-int}). (The minimum exists since $\rho$ and therefore $\sigma$ vary continuously as the the left-hand side of (\ref{imp-def}) varies.)

Note also that since $g$ is continuous it is bounded on $[\rho_{\rm min},\rho_{\rm max}]$, and that $n-k-t$ grows linearly with $n$.
It follows that for all sufficiently large $n$ (depending on $p, q, M$) we have
\begin{equation} \label{inq-largen}
\frac{1 \pm \frac{g(\rho)}{n-k-t}}{\sigma \sqrt{2\pi(n-k-t)}} \leq \frac{c}{\sqrt{n}}.
\end{equation}
where $c=c(p,q,M)$ is a constant. 
Since for fixed $M$ there are only finitely options for $p, q$, we can find a single constant $c=c(M)$ so that for all large enough $n$ (\ref{inq-largen}) holds for all $k, t$ under consideration.  

Combining all of these bounds, it follows that for all large enough $n$ and for all $k, t$ satisfying $k \leq 0.49n$ and $0 \leq t \leq k+1$ we have
$$
\frac{f(n,k,t)}{\binom{n-k+1}{k}} 
\leq \frac{cM(n+1)}{(M-p)\sqrt{n}}\left(\frac{ABCDEIJ}{FG}\right)^n
$$
where $p, q$ are associated with $k$ and $t$ as described earlier, and $A, B$, et cetera are the various terms (depending on $p, q$ and $M$) that we have just defined. It follows that
$$
\sum_{t=0}^{k+1}\frac{f(n,k,t)}{\binom{n-k+1}{k}} 
\leq \frac{cM(n+1)}{(M-p)\sqrt{n}}\left(\max \left\{\frac{ABCDEIJ}{FG}:0 \leq q \leq p+1\right\}\right)^n
$$
(taking $q$ up to $p+1$ in the max calculation is necessary since $t$ ranges up to $k+1$). If there is an $M$ such that
\begin{equation} \label{inq-comp}
\max \left\{\frac{ABCDEIJ}{FG}:0 \leq q \leq p+1\right\} < 1
\end{equation}
for all $0.1M \leq p \leq \weakbestinc M$, then 
we obtain (\ref{eq-inc-goal}) (equivalently (\ref{main-inq})), that is,
$$
\sum_{t=0}^{k+1} \frac{f(n,k,t)}{\binom{n-k+1}{k}}
\leq \frac{1}{n^2\log n} 
$$
for all sufficiently large $n$, completing the proof of Theorem \ref{thm-random-tree-inc-head}.

To make the computation manageable, we proceed in stages. We can begin, for example, by showing that with $M=100$, (\ref{inq-comp}) holds for $10 \leq p \leq 23$, yielding that the independent set sequence of the uniform labelled tree on $[n]$ is a.a.s. increasing up to $0.23n$. So from here on we may restrict attention to $p \geq 0.23M$. With $M=1000$, (\ref{inq-comp}) holds for $230 \leq p \leq 274$, allowing us in the sequel to restrict to $p \geq 0.274M$. Bootstrapping in this way, we eventually get to $M=7500$, at which value (\ref{inq-comp}) holds for $p \leq 2100$, yielding the bound claimed in Theorem \ref{thm-random-tree-inc-head}. All computations were performed on {\tt Mathematica}.

Now we turn to Theorem \ref{thm-random-tree-dec-tail}. Here our goal is to establish
\begin{equation} \label{dec-goal}
\sum_{t=k+1}^{n-k} \frac{f(n,k,t)}{\binom{n-k}{k+1}}
\leq \frac{1}{n^2\log n} 
\end{equation}
for all large enough $n$, and $k \geq \weakbestdec n$. By the discussion immediately preceding the statement of Theorem \ref{thm-random-tree-dec-tail} (specifically, by combining (\ref{meir-moon}) and Theorem \ref{thm-decreasing-tail}), we may assume $k \leq 0.362n$.

We need to make two minor changes to the simple bounds established in the proof of Theorem \ref{thm-random-tree-inc-head}:
\begin{itemize}
\item In bounding $\binom{n-k}{t} \leq B(p,q,M)^n$ we used $t/(n-k) \leq q/(M-p) \leq 1/2$
which allowed the conclusion $H(t/(n-k) \leq H(q/(M-p)))$. Now we no longer have $t/(n-k) \leq 1/2$, so instead we say
$$
\frac{q-1}{M-(p-1)} \leq \frac{t}{n-k} \leq \frac{q}{M-p}.
$$
If $1/2 \not \in [(q-1)/(M-(p-1)),q/(M-p)]$ then 
$$
H\left(\frac{t}{n-k}\right) \leq \max\left\{H\left(\frac{q-1}{M-(p-1)}\right), H\left(\frac{q}{M-p}\right)\right\},
$$
while if $1/2 \in [(q-1)/(M-(p-1)),q/(M-p)]$ then $H(t/(n-k)) \leq 1$. So we may replace $B(p,q,M)$ with
$$
B'(p,q,M) = \exp_2\left\{\left(1-\frac{(p-1)}{M}\right)b(p,q,M)\right\}
$$
where
$$
b(p,q,M) = \left\{
\begin{array}{ll}
\max\left\{H\left(\frac{q-1}{M-(p-1)}\right), H\left(\frac{q}{M-p}\right)\right\} & \mbox{if $\frac{1}{2} \not \in \left[\frac{q-1}{(M-(p-1)},\frac{q}{M-p}\right]$} \\
1 & \mbox{otherwise.}
\end{array}
\right.
$$

\item Because the denominator on the left-hand side of (\ref{dec-goal}) is $\binom{n-k}{k+1}$ rather than $\binom{n-k+1}{k}$, we replace $G(p,M)$ with
$$
G'(p,M) = \exp_2\left\{\left(1-\frac{p}{M}\right)H\left(\frac{p+1}{M-p}\right)\right\}.
$$
\end{itemize}

We also need to deal with an issue in the application of Theorem \ref{thm-good}, which required that $k/(n-k-t)$ be bounded away from infinity. This was the case when $k\leq \weakbestinc$ and $t \leq k+1$, but it is no longer the case when $k \geq \weakbestdec$ and $t \geq k+1$, since $t$ (the possible number of extensions of an independent set of size $k$ to one of size $k+1$) can be as large as $n-k-1$.  

One way to get around this problem is to use an alternate, slightly weaker, upper bound on  
$$
\sum_{j=1}^k {k \brace j}(n-k-t)_j \left(\frac{n}{k}\right)^j
$$
when $t$ is close to $n-k$. Using the standard identity
$\sum_{\ell=1}^N {N \brace \ell}(x)_\ell = x^N$ we have
$$
\frac{1}{n^k}\sum_{j=1}^k {k \brace j}(n-k-t)_j\left(\frac{n}{k}\right)^j \leq \left(\frac{n-k-t}{k}\right)^k \leq \left(\frac{M-(p+q-2)}{p-1}\right)^\frac{(p-1)n}{M}.
$$
It follows that for $t$ close to $n-k$ we may replace $C(IJ)^n/\sqrt{n}$ with $K^n$ where
$$
K(p,q,M) = \left(\frac{M-(p+q-2)}{p-1}\right)^\frac{p-1}{M}.
$$
Specifically, we may bound
$$
\frac{f(n,k,t)}{\binom{n-k}{k+1}}  
\leq \left\{
\begin{array}{ll}
\frac{CM(n+1)}{(M-p)\sqrt{n}}\left(\frac{AB'CDEIJ}{FG'}\right)^n & \mbox{if $t \leq 0.99n-k$} \\
\frac{M(n+1)}{(M-p)}\left(\frac{AB'CDEK}{FG'}\right)^n & \mbox{if $t > 0.99n-k$} 
\end{array}.
\right.
$$

The computational verification of (\ref{dec-goal}) now proceeds in a very similar manner to that of (\ref{eq-inc-goal}), and we omit the details.

\subsubsection{Deriving (\ref{eq-enkt})} \label{sec-mtt}

Here we use the Matrix Tree Theorem to find an explicit expression for $e(n,k,t)$, the probability that, in a uniformly chosen labelled tree on $[n]$, a particular set of size $k$ is independent and has exactly $t$ extensions to an independent set of size $k+1$.

Given two disjoint subsets $K, L$ of $[n]$ with $|K|=k \geq 1$ and $|L|=\ell$, denote by $T_{K,L}$ the set of trees on $[n]$ with $K$ an independent set and with $L$ having no edges to $K$.

\begin{claim} \label{clm-mtt}
$$
|T_{K,L}| = n^{n-2}\left[\left(1-\frac{k}{n}\right)^{\ell-1}\left(1-\frac{(k+\ell)}{n}\right)^k\right].
$$
\end{claim}

\medskip

\noindent {\bf Proof}: $T_{K,L}$ is exactly the set of spanning trees of the graph $G(K,L)$ obtained from $K_{[n]}$, the complete graph on vertex set $[n]$, by deleting all the edges inside $K$, as well as all edges from $L$ to $K$. 

The Laplacian of $G(K,L)$, with the rows and columns indexed first by vertices in $K$, then $L$, then the rest of the vertices (call this set $M$), is a block matrix. 
\begin{itemize}
    \item The block with rows indexed by $K$, columns indexed by $K$, has $0$'s off the diagonal, and $n-k-l$'s down the diagonal.
    \item The block with rows indexed by $K$, columns indexed by $L$, is all $0$.
    \item The block with rows indexed by $K$, columns indexed by $M$, is all $-1$.
    \item The block with rows indexed by $L$, columns indexed by $L$, has $-1$'s off the diagonal, and $n-k-1$'s down the diagonal.
    \item The block with rows indexed by $L$, columns indexed by $M$, is all $-1$.
    \item The block with rows indexed by $M$, columns indexed by $M$, has $-1$'s off the diagonal, and $n-1$'s down the diagonal.
\end{itemize}
(No other blocks need be specified --- the matrix is symmetric). This matrix has
\begin{itemize}
    \item $0$ as an eigenvalue with geometric multiplicity at least $1$ (all row sums are $0$);
    \item $n-k-\ell$ as an eigenvalue with geometric multiplicity at least $k$ (on subtracting $n-k-\ell$ from each diagonal entry, the first $k$ rows become identical, and the sum of the rows indexed by $L$ is a multiple of this common value);
    \item $n-k$ as an eigenvalue with geometric multiplicity at least $\ell-1$ (on subtracting $n-k$ from each diagonal entry, the $\ell$ rows indexed by $L$ become identical); and
    \item $n$ as an eigenvalue with geometric multiplicity at least $n-k-\ell$ (on subtracting $n$ from each diagonal entry, the $n-k-\ell-1$ rows indexed by $M$ become identical, and the sum of the remaining rows is a multiple of this common value).
\end{itemize}
Since $1+k+(\ell-1)+(n-k-\ell)=n$ it follows that these lower bounds on geometric multiplicities are equalities, and that the algebraic multiplicities of all the eigenvalues coincide  with their geometric multiplicities. So from the Matrix Tree Theorem we get
\begin{eqnarray*}
|T_{K,L}| & = & n^{n-k-\ell-1}(n-k)^{\ell-1}(n-k-\ell)^k \\
& = & n^{n-2}\left[\left(1-\frac{k}{n}\right)^{\ell-1}\left(1-\frac{(k+\ell)}{n}\right)^k\right].
\end{eqnarray*}
\qed

\medskip

Now let $\emptyset \neq K \subseteq [n]$ be given, as well as $T \subseteq [n]\setminus K$ ($T$ might be empty). Set $|K|=k$ and $|T|=t$.

\begin{claim} \label{clm-mtt2}
The number of trees on $[n]$ with $K$ as an independent set, and with $T$ as the {\em exact} set of vertices that extend $K$ to an independent set of size $k+1$, is
$$
n^{n-2}\sum_{\ell=0}^{n-k-t} (-1)^\ell \binom{n-k-t}{\ell}\left(1-\frac{k}{n}\right)^{\ell+t-1}\left(1-\frac{(k+t+\ell)}{n}\right)^k.
$$
\end{claim}

\noindent {\bf Proof}: Let $U_{K,T}$ be the set of trees with $K$ as an independent set and with $T$ {\em among} the set of vertices that extend $K$ to an independent set of size $k+1$; we know from Claim \ref{clm-mtt} that
$$ 
|U_{K,T}| = n^{n-2}\left(1-\frac{k}{n}\right)^{t-1}\left(1-\frac{(k+t)}{n}\right)^k.
$$ 
Let the vertices of $[n]\setminus (K \cup T)$ be $v_1, \ldots, v_{n-k-t}$. Let $A_j$ be the set of trees in $U_{K,T}$ in which there is no edge from $v_j$ to $K$. Then the number of trees on $[n]$ with $K$ as an independent set, and with $T$ as the {\em exact} set of vertices that extend $K$ to an independent set of size $k+1$, is
$$
U_{K,T} \setminus (A_1 \cup \cdots \cup A_{n-k-t}). 
$$
If $L$ is any subset of $\{1, \ldots, n-k-t\}$ then $\cap_{i \in L} A_i$ is exactly the set of trees with $K$ independent, and with $T \cup L$ among the set of vertices that extend $K$ to an independent set of size $k+1$, so by Claim \ref{clm-mtt} we have
$$
\left|\cap_{i \in L} A_i\right| = n^{n-2}\left(1-\frac{k}{n}\right)^{t+\ell-1}\left(1-\frac{(k+t+\ell)}{n}\right)^k.
$$
So, by inclusion-exclusion, the number of trees on $[n]$ with $K$ as an independent set, and with $T$ as the {\em exact} set of vertices that extend $K$ to an independent set of size $k+1$, is
\begin{eqnarray*}
 & n^{n-2}\left(1-\frac{k}{n}\right)^{t-1}\left(1-\frac{(k+t)}{n}\right)^k & \\
 & - n^{n-2} \sum_{\ell=1}^{n-k-t} (-1)^{\ell-1} \binom{n-k-t}{\ell}\left(1-\frac{k}{n}\right)^{t+\ell-1}\left(1-\frac{(k+t+\ell)}{n}\right)^k & 
\end{eqnarray*}
or more compactly
$$
n^{n-2} \sum_{\ell=0}^{n-k-t} (-1)^\ell \binom{n-k-t}{\ell}\left(1-\frac{k}{n}\right)^{t+\ell-1}\left(1-\frac{(k+t+\ell)}{n}\right)^k.  
$$
\qed

The claimed expression (\ref{eq-enkt}) for $e(n,k,t)$ (the probability that in a uniformly chosen labelled tree on $[n]$ a given set of size $k$ is independent and has exactly $t$ extensions to an independent set of size $k+1$) follows from Claim \ref{clm-mtt2} by first summing over all possible choices for $T$ (the set of extensions) and then using Cayley's formula.

\medskip

\noindent {\bf Acknowledgement}: We are grateful to the referee who provided us with the combinatorial proof of (\ref{id-Boyadzhiev}), and who also made many helpful suggestions to improve the clarity of our presentation.

\end{document}